\documentclass {article}

\usepackage{amsfonts,amsmath,latexsym}
\usepackage{amstext, amssymb, amsthm}
\usepackage[scanall]{psfrag}
\usepackage {graphicx}
\usepackage{verbatim}

\newcommand {\R}{\mathbb{R}}

\newcommand {\D}{\mathbb{D}}

\newcommand {\del}{\partial}

\newcommand {\Area}{\operatorname{Area}}
\newcommand {\vol}{\operatorname{Vol}}

\newcommand {\calC}{\mathcal{C}}
\newcommand {\calR}{\mathcal{R}}

\newtheorem {thm} {Theorem}
\newtheorem* {thm*} {Theorem}
\newtheorem {prop} [thm] {Proposition}
\newtheorem {lemma} [thm] {Lemma}
\newtheorem {cor} [thm] {Corollary}
\newtheorem {defn} {Definition}
\newtheorem {rmk} {Remark}

\begin {document}

\title{ 
{\bf Interpolating Between Torsional
Rigidity and Principal Frequency}}
\author{Tom Carroll and Jesse Ratzkin \\
University College Cork \\
{\tt t.carroll@ucc.ie} and {\tt j.ratzkin@ucc.ie}} 
\maketitle

\begin {abstract} 
\noindent A one-parameter family of 
variational problems is introduced that interpolates between 
torsional rigidity and the first Dirichlet eigenvalue of the Laplacian.
The associated partial differential equation is derived, which is shown
to have positive solutions in many cases. 
Results are obtained regarding extremal domains and regarding 
variations of the domain or the parameter. 
\end {abstract}

\section{Introduction}

Various geometric and physical constants can be associated with 
a bounded domain in the plane: perimeter, area, diameter, 
transfinite diameter or capacity, torsional rigidity and principal frequency. 
We focus on the latter two physical quantities. 
The torsional rigidity $P(D)$ of a bounded domain $D$ in the plane is given by 
\begin {equation} \label{variational-torsion} 
\frac{4}{P(D)} = \inf \left \{ \frac{ \int_D |\nabla u|^2 
dA}{\left ( \int_D u dA \right )^2 }: u \not
\equiv 0, \left. u \right|_{\del D} = 0 \right \}
\end {equation}
and the principal frequency $\lambda(D)$ is given by 
\begin {equation} \label {variational-frequency}
\lambda(D) = \inf \left \{ \frac{ \int_D |\nabla u|^2 
dA}{ \int_D u^2 dA  }: u \not \equiv 0, 
\left. u \right|_{\del D} = 0 \right \}.
\end {equation} 
Each has a physical significance. In case $D$ is simply connected, the torsional 
rigidity is a measure of the torque required per unit length per unit angle of twist 
when a beam with cross section $D$ undergoes torsion.  
The principal frequency is a measure of the lowest note that a drum 
of shape $D$ can produce. An extensive literature is devoted to each.  
The classic book \cite{PS} by P\' olya and Szeg\H o contains a wealth of information on
both, in particular the variational characterizations \eqref{variational-torsion} 
and \eqref{variational-frequency} and the Euler-Lagrange partial 
differential equations (PDEs) listed below. 
The corresponding Euler-Lagrange equations 
for these two problems are 
\begin {equation} \label {torsion-pde}
\Delta u + 2 = 0, \qquad \left. u 
\right|_{\del D} = 0
\end {equation}
for torsional rigidity and 
\begin {equation} \label {frequency-pde}
\Delta u + \lambda u = 0, \qquad \left. u 
\right|_{\del D} = 0
\end {equation}
for frequency. Classical theorems tell us that both these equations 
have positive solutions for a bounded domain $D$, provided 
the boundary $\del D$ has some regularity. 

One can measure the fundamental frequency of a bounded domain in any dimension;
both the variational formulation (\ref{variational-frequency}) and the PDE
(\ref{frequency-pde}) remain the same. 
The torsion function, the solution to the problem \eqref{torsion-pde}, is
of probabilistic significance in any dimension, being the expected 
exit time of Brownian motion from the domain $D$ relative to Wiener measure. 
 
The purpose of this note is to introduce a one-parameter family of variational 
problems and corresponding PDEs that interpolate between the torsional 
rigidity and principal frequency problems listed above, and go somewhat beyond 
them. 
In the remainder of this introduction, we briefly describe the interpolating problem and 
the relevant physical quantities, state some results, and outline the remainder 
of the paper.  

Let $D\subset \R^n$, $n\geq 2$, be a bounded domain with smooth boundary 
({\it i.e.\/} one can locally write the boundary of $D$ as the graph of a $C^\infty$ function). 
\begin {defn} \label{variation-p}
For each $p\geq 1$, we define $\calC_p(D)$ by 
\begin {equation} \label{C-p}
\calC_p(D) = \inf \left \{ \Phi_p(u) = 
\frac{\int_D |\nabla u|^2 dV}{\left ( 
\int_D u^p dV\right )^{2/p} } : u \in L^p(D) 
\cap W^{1,2}_0(D), u \not \equiv 0 \right \}.
\end {equation}
\end {defn}  
\noindent In particular, $4 / \calC_1(D)$ is the torsional rigidity and 
$\calC_2(D)$ is the fundamental frequency. 
The scaling law for $\calC_p$ is 
\begin{equation}\label{scaling-law}
\calC_p(rD) = r^{n - 2 -\frac{2n}{p}} \calC_p(D), \quad r>0.
\end{equation}
In Section~\ref{euler-eqn-sec} we derive the 
corresponding Euler-Lagrange equation. 
%
% Theorem 1:  Euler-Lagrange equation 
%
\begin {thm}\label{thm1}
Let  $D\subset \R^n$ be a bounded domain with smooth boundary. 
Let $p \geq 1$. The critical points of the functional 
$$
\Phi_p(u) = \frac{\int_D |\nabla u|^2 dV}
{\left ( \int_D u^p dV \right )^{2/p}}
$$
in $\displaystyle  L^p(D) \cap W^{1,2}_0(D)$ satisfy the 
PDE  
\begin {equation} \label {p-pde} 
\Delta u + \Lambda u^{p-1} = 0, 
\qquad \left. u \right|_{\del D} = 0
\end {equation} 
for some constant $\Lambda$. 
\end {thm} 
After verifying that this is the correct Euler-Lagrange equation, 
we discuss solvability of the PDE, borrowing some classical results 
of Pohozaev \cite{Poh}, and some basic properties of solutions. 

In contrast to the functional $\Phi_p$, the differential equation (\ref{p-pde}) 
is not scale-invariant unless $p=2$. 
If $u$ solves (\ref{p-pde}) for some constant $\Lambda$ 
and $k>0$ is a constant, then $v=ku$ satisfies 
$$
\Delta v + k^{2-p} \Lambda\, v^{p-1} = 0.
$$
Conversely, given a solution $u$ to (\ref{p-pde}) and a constant $\alpha>0$, 
we see that 
$$
v = \left ( \frac{\alpha}{\Lambda} \right 
)^{\frac{1}{2-p}} u \qquad \text {solves} 
\qquad \Delta v + \alpha\, v^{p-1} = 0.
$$
Thus, if $p\neq 2$, by rescaling we obtain solutions to the equation $\Delta u + \Lambda u^{p-1} = 0$
for all $\Lambda >0$. 
The variational problem, however, will select a particular Lagrange multiplier $\Lambda$. 

\begin {lemma}\label{lemma2}
Let $\calC_p(D)$ be given by  \eqref{C-p} and let $u$ be a positive solution to 
\eqref{p-pde} for some $\Lambda >0$. Then 
\begin {equation} \label{p-eigenval-lagrange}
\calC_p(D) = \Lambda 
\left ( \int_D u^p dV \right )^{(p-2)/p}.
\end {equation} 
\end {lemma} 
Notice that the right hand side of (\ref{p-eigenval-lagrange}) is invariant
under scaling of the function $u$. 
We recover a classical formula for the torsion function 
when $p=1$, namely $P(D) = 2\int_D u\,dA$ where $u$ is the solution of 
\eqref{torsion-pde}. The case $p=2$ returns the first Dirichlet 
eigenfunction of the Laplacian. 

Motivated by this lemma and torsional 
rigidity, we make the following definition. 
% definition of p-torsional rigidity
\begin {defn} Let $D\subset \R^n$ be a bounded smooth domain and either $n=2$, 
or $n\geq 3$ and $p<\frac{2n}{n-2}$.  Then the $p$-torsional rigidity is given by 
\begin {equation} \label{p-torsion}
\calR_p(D) = \frac{4}{\Lambda} \left (
\int_D u^p dV \right )^{(2-p)/p},
\end {equation}
where $u>0$ solves the boundary value 
problem (\ref{p-pde}).
\end {defn} 
\noindent Then $\calR_1(D) = P(D)$ and $\calR_2(D) = 4/\lambda(D)$. We see immediately from the definition 
and (\ref{p-eigenval-lagrange}) that
\begin {equation} \label {p-lambda-eqn}
\calC_p(D)\, \calR_p(D) = 1, \quad p\geq 1.
\end {equation}

Next, in Section~\ref{examples-sec}, we examine the specific case where $D$ is
a round ball or an infinite slab $\{-1 < x_n 
< 1\} \subset \R^n$. 
We prove some comparison results for 
$\calC_p(D)$ in Section \ref{comparison-sec}, 
varying either the domain 
$D$ or $p$. In particular, we prove the 
following theorem. 
\begin {thm} \label{holder-comparison}
Let $D \subset \R^n$ be a bounded domain 
with volume $V = \vol(D)$, and let $1 \leq 
p < q$. Then 
\begin{equation}\label{FK}
V^{2/p} \calC_p(D) > V^{2/q} \calC_q (D).
\end{equation}
\end {thm} 
\noindent In dimension 2, the case $p=1$ and $q=2$ of the inequality 
\eqref{FK} becomes $\lambda(D) < 4A(D)/ P(D)$ which one may find on 
Page 91 of \cite{PS} and  relates the fundamental frequency of a domain 
to its area and its torsional rigidity. 
The scaling law \eqref{scaling-law} shows that 
\begin{equation}\label{scaling-law2}
V(rD)^{2/p} \calC_p(rD) = r^{n-2} V(D)^{2/p} \calC_p(D), \quad p\geq 1,
\end{equation}
which agrees with the classical scaling laws for the torsional rigidity and the 
principal frequency. 
In Section \ref{extreme-sec} we characterize extremal domains for 
$\calC_p$ under certain conditions. 
We prove an inequality of  Faber-Krahn type. See \cite{F,K} for the classical 
Faber-Krahn inequality for principal frequency (the case  $p=2$) and see 
\cite{P} or Appendix A of \cite{PS} for P\'olya's proof of the Saint-Venant theorem that among all simply connected domains of given area the disk has 
the largest torsional rigidity. 
%
% our Faber-Krahn theorem
%
\begin {thm} \label {faber-krahn}
Let $p \geq 1$, let $D \subset \R^n$ be a bounded 
domain with $C^\infty$ boundary, and 
let $B\subset \R^n$ be the round 
ball, centered at the origin, with the 
same volume as $D$. Then 
$\calC_p(D) \geq \calC_p(B)$. Moreover, equality 
only occurs if $D=B$ almost everywhere. 
\end {thm} 
\noindent We also discuss convex domains of fixed inradius in this section. 
We conclude with a short list of open question in Section~\ref{open-sec}.

\section {The variational problem and its 
Euler-Lagrange equation} 
\label {euler-eqn-sec} 

We take $p\geq 1$ and a bounded domain $D
\subset \R^n$ with a $C^\infty$ boundary,
and for $u \in L^p(D) \cap 
W^{1,2}_0(D)$ not identically zero 
define the functional 
$$\Phi_p(u) = \frac{\int_D |\nabla u|^2 
dV}{\left ( \int_D u^p dV \right )^{2/p}} 
= \frac{ \|\nabla u\|_{L^2(D)}^2} {\| u
\|_{L^p(D)}^2}.$$ Our first task is 
to derive the Euler-Lagrange equation 
(\ref{p-pde}). 

\begin {proof} Observe that $\Phi_p$ is scale-invariant; that is, if $k>0$ then 
$\Phi_p(ku) = \Phi_p(u)$. Thus, we can reformulate the condition that $u$ 
is a critical point of $\Phi_p$ as a constrained critical point problem: find 
the critical points of $\int_D |\nabla u|^2 dV$ subject to the contraint 
$\int_D u^p dV = 1$. 
Any constrained critical point must satisfy 
$$
\left. \frac{d}{d\epsilon} \right|_{\epsilon 
= 0} \int_D |\nabla (u+\epsilon v)|^2 dV 
= \Lambda \left. \frac{d}{d\epsilon} 
\right |_{\epsilon= 0} \int_D (u+\epsilon 
v)^p dV
$$
for all $v \in L^p(D) \cap W^{1,2}_0(D)$, 
where $\Lambda$ is the Lagrange multiplier. 
Next recall that $C^\infty_0(D)$ 
is dense in both $L^p(D)$ and $W^{1,2}_0(D)$
(see, for example, Section~7.6 of \cite{GT}), 
so without loss of generality we can take 
$u$, $v \in C^\infty_0(D)$. Thus we can freely 
differentiate underneath the integral sign, 
and a quick computation shows that the equation above is 
equivalent to 
\begin {equation} \label{integral-euler-eqn}
0 = \int_D -\langle \nabla u, \nabla v\rangle 
+ \Lambda u^{p-1} v dV = 
\int_D v(\Delta u + \Lambda u^{p-1}) dV.
\end {equation}
Here we have absorbed a factor of $2$ and 
a factor of $p$ into the Lagrange multiplier 
$\Lambda$. If equation (\ref{integral-euler-eqn})
is to hold for all compactly supported $v$ in 
$D$, then $u$ must satisfy the PDE
$$\Delta u + \Lambda u^{p-1} = 0$$
for some constant $\Lambda$ as claimed. 
\end {proof} 

\begin {rmk} This is a familiar differential 
equation, occurring (for instance) in 
the study of scalar curvature under a conformal 
change of metric; see \cite{LP}. 
\end {rmk} 

\begin {rmk} One can equally well study the 
functional 
$$u \mapsto \frac{\left ( \int_D |\nabla u|^q 
dV \right )^{2/q}}{\left (\int_D u^p dV 
\right )^{2/p}}.$$
In this case, the Euler-Lagrange equation is 
$$0 = \Delta_q u + \Lambda u^{p-1} = 
\operatorname{div} (|\nabla u|^{q-2} \nabla u)
+ \Lambda u^{p-1}.$$
This differential equation is either
singular (for $q<2$) or degenerate (for 
$q>2$) at the critical points of $u$. 
Thus, we do not expect to 
have as well-developed a theory attached 
to the more general variational problem.
\end {rmk} 

To examine the solvability of equation 
(\ref{p-pde}), we first recall the following  
classical theorem of Pohozaev \cite{Poh}.
\begin {thm*} (Pohozaev) Let $D\subset \R^n$ be 
a bounded domain with smooth boundary, and let $f(u)$ 
be a Lipschitz function. If $n=2$ and $f$ satisfies 
the estimate 
$$|f(u)| \leq A + B |u| e^{c|u|^a}, \qquad a<2,$$
then one can find eigenfunctions of the PDE 
$$ \Delta u + \lambda f(u) = 0, \qquad \left. 
u \right|_{\del D} = 0.$$
If $n\geq 3$ and $f$ satisfies the estimate
$$|f(u)| \leq A + B|u|^m, \qquad m < 
\frac{n+2}{n-2},$$ 
then again one can find eigenfunctions of 
the PDE 
$$\Delta u + \lambda f(u) = 0, \qquad 
\left. u \right|_{\del D} = 0.$$
Conversely, if $D$ is star-shaped with respect 
to the origin and $u \geq 0$, not identically 
zero, solves 
$$\Delta u + u^m = 0, \qquad \left. 
u \right|_{\del D}  = 0$$
then $m < \frac{n+2}{n-2}$. 
\end {thm*} 

Applying Pohozaev's theorem, we immediately obtain the following corollary. 
\begin {cor} 
There is a positive solution to the Euler-Lagrange equation (\ref{p-pde})
if either $n=2$, or $n\geq 3$ and $p < \frac{2n}{n-2}$. 
On the other hand, if $n\geq 3$, $p> \frac{2n}{n-2}$, and 
$D \subset \R^n$ is star-shaped,  then equation (\ref{p-pde}) does 
not have a positive solution. 
\end {cor} 
It is well-known that for the critical value $p=\frac{2n}{n-2}$ 
a minimizing sequence for the function $\Phi_p$ 
will typically become unbounded. Thus it 
is difficult to determine whether one can 
realize the corresponding infimum $\calC_p(D)$ as $\Phi_p(u)$ for 
some $u >0$. One can find a treatment of this blow-up 
phenomenon, which reflects the loss of 
compactness in the Sobolev embedding theorem,
in \cite{Tr} and Section~4 of \cite{LP}. 

We use the maximum principle to prove 
the following lemma. 
\begin {lemma} Let $D \subset \R^n$ be 
a smooth, bounded domain and $\Lambda$ 
a constant. Then there is at most one 
positive solution to the boundary value 
problem (\ref{p-pde}). 
\end {lemma} 

\begin {proof} Suppose $u$ and $v$ are 
distinct positive solutions to (\ref{p-pde}), 
and let $D^*$ be a connected component of 
$\{u \neq v\}$. Without loss of generality
we can assume that $\Lambda >0$ and that 
$u>v$ on $D^*$. We can also assume $\del D^*$ 
is smooth by Sard's theorem. Then $u-v$ 
is a positive solution to the 
boundary value 
$$\Delta (u-v) = \Lambda (u^{p-1} - v^{p-1}) 
> 0, \qquad \left. u \right|_{\del D^*} 
= \left. v \right|_{\del D^*},$$
which contradicts the maximum principle. 
\end {proof} 

If $D$ is a convex domain, we can in fact say more. A theorem of Korevaar (see 
Theorem 2.5 of \cite{Kor} and the remark immediately following it) implies
\begin {cor} 
If $D\subset \R^n$ is a bounded, smooth, strictly convex domain and $u > 0$ solves the boundary value problem (\ref{p-pde}) then $v = -\log(u)$ is convex. 
\end {cor} 

We complete this section by proving Lemma~\ref{lemma2}. 
\begin {proof}[Proof of Lemma~\ref{lemma2}]By Theorem~\ref{thm1}, 
$\calC_p(D) = \Phi(u)$. We integrate by parts and use \eqref{p-pde}: 
\begin {eqnarray*} 
\int_D |\nabla u|^2 dV & = & \int_D \langle \nabla u, \nabla u\rangle dV 
= -\int_D u\Delta u dV \\
& = & \Lambda \int_D u^p dV 
= \Lambda \left ( \int_D u^p dV \right )^{2/p}
\left ( \int_D u^p dV \right )^{(p-2)/p}.
\end {eqnarray*}
Rearranging yields  
\begin{equation*}
\calC_p(D) = \frac{\int_D |\nabla u|^2dV}{\left(\int_D u^p\right)^{2/p}} 
= \Lambda \left ( \int_D u^p dV \right )^{(p-2)/p}.
\tag*{\qedhere}
\end{equation*}
\end {proof} 

\section {Examples}
\label {examples-sec}

We begin with the case of an infinite slab
$$
S = \{ (x_1, \dots, x_n) : -1 < x_n < 1 \}.
$$ 
In order that the variational problem makes sense, one can truncate to 
obtain $S$ as a limit of $D_R$ as $R \rightarrow \infty$, where 
$$
D_R = \{ (x_1, \dots, x_n) : -1 < x_n < 1, -R < x_j < R\},
$$
and in the limit we recover the same Euler-Lagrange 
equation $\Delta u + \Lambda u^{p-1} = 0$.
We look for a solution which depends only on $x_n$, which will solve 
the following boundary value problem for an ordinary differential 
equation (ODE):
\begin {equation} \label {slab-ode} 
u'' + \Lambda u^{p-1} = 0 , \qquad 
u(-1) = 0 = u(1).
\end {equation} 

A quick computation shows 
$$\frac{d}{dt} \left ( (u')^2 + \frac{2}{p} 
\Lambda u^p \right ) = 2u' (u'' + \Lambda 
u^{p-1}) = 0,$$
so in phase space solutions to (\ref{slab-ode})
will lie on level sets of the energy function 
\begin {equation} \label {slab-1st-integral} 
E = (u')^2 + \frac{2\Lambda}{p} u^p.
\end {equation} 
Equivalently, for any solution $u$ to 
(\ref{slab-ode}) there is a constant $E$ 
such that 
$$u' = \sqrt{E - \frac{2\Lambda}{p} u^p}.$$
One can use this last equation to write 
down all the solutions to (\ref{slab-ode}) 
up to quadrature, or in terms of hypergeometric 
functions. We sketch some level sets of the 
energy function (\ref{slab-1st-integral}) below. 

\begin {figure}[h]
\begin {center}
\includegraphics [width=4in]{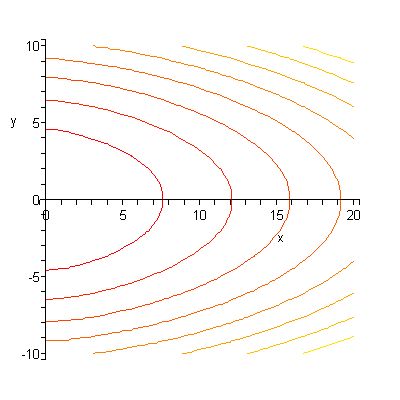}
\label{slab-level-sets}
\caption{Level sets for the energy function
(\ref{slab-1st-integral}) with $p=3/2$.}
\end {center}
\end {figure} 

Indeed, we can compute $\calC_p(S)$ for the slab using just the knowledge 
that a positive solution to (\ref{p-pde}) depends only on $x_n$. 
On the truncated domain $D_R$, we have 
$$
\int_{D_R} |\nabla u|^2 dV = (2R)^{n-1} 
\int_{-1}^1 (\del_{x_n} u)^2 dx_n = C_1 R^{n-1},
$$
while 
$$
\left ( \int_{D_R} u^p dV \right )^{2/p} = (2R)^{\frac{2}{p}(n-1)}
\left ( \int_{-1}^1 u^p(x_n) dx_n \right )^{2/p} = C_2 R^{\frac{2}{p}(n-1)}.
$$
Taking a ratio we see that $\Phi_p(u) = \mathcal{O} (R^{(n-1)(1- \frac{2}{p})})$, 
and so
$$
\calC_p(S) = \left \{ \begin {array}{cc} 0 & 1 \leq p < 2 \\
\frac{\pi^2}{4} & p=2\\
\infty & p>2, \end {array}\right. 
$$
where in the $p=2$ case we have listed the (well-known) value of the first
Dirichlet eigenvalue of the Laplacian of a slab of width $2$. 

\begin {rmk} 
Notice that the solution for a slab does not depend on dimension at all. 
\end {rmk} 

Next we examine the case of the ball. By a moving planes argument \cite{GNN},
any positive solution to (\ref{p-pde}) is radial. 

We treat the cases $n=2$ 
and $n\geq 3$ separately. First consider 
the unit disc in the plane. A radial 
solution $u$ will satisfy the ODE 
$$u'' + \frac{1}{r} u' + \Lambda u^{p-1} 
= 0, \qquad u(1) = 0.$$
Change variables by $t = -\log (r)$ and 
let $v(t) = u(e^{-t})$ to get the boundary 
value problem 
\begin {equation} \label {disc-ode} 
e^{2t} \ddot{v} + \Lambda v^{p-1} = 0, 
\qquad v(0) = 0,\end {equation} 
where the $\lim_{t \rightarrow \infty} v(t)$ 
exists and is a positive number. Here we 
have used a dot to denote differentiation 
with respect to $t$ and a dash to denote 
differentiation with respect to $r$. 

In dimension $n \geq 3$ a radial solution will 
satisfy the ODE 
$$u'' + \frac{n-1}{r} u' + \Lambda u^{p-1}  
= 0, \qquad u(1) = 0.$$
We do a similar change of variables, but this 
time rescale $v$ by a radial factor, letting 
$u = e^{-kt} v(t)$, where again $t = -\log (r)$. 
Again, we will use a dot to denote differentiation 
with respect to $t$ and a dash to denote 
differentiation with respect to $r$. Under this 
change of variables, the ODE becomes 
\begin {equation} \label{ball-ode-a} 
e^{(2-k)t} (\ddot v + (n+2k-2)\dot v + 
k(n+k-2) v) + \Lambda e^{-k(p-1)t} v^{p-1} 
= 0, \qquad v(0) = 0,\end {equation}
where $\lim_{t \rightarrow \infty} e^{-kt} v(t)$
exists and is a positive number. We can 
eliminate the $\dot v$ term by choosing $k=
\frac{2-n}{2}$, and so (\ref{ball-ode-a}) is 
now 
\begin {equation} \label {ball-ode-b}
e^{\frac{(n+2)t}{2}} \left ( \ddot v  - 
\left (\frac{n-2}{2} \right )^2 v \right ) 
+ \Lambda e^{\frac{(n-2)(p-1)t}{2}} v^{p-1}  
= 0. 
\end {equation} 
In the critical case of $p = \frac{2n}{n-2}$, the exponential terms coincide, 
and we obtain the familiar differential equation 
$$ 
\ddot v - \left ( \frac{n-2}{2} \right )^2 v + \Lambda v^{\frac{n+2}{n-2}} = 0,
$$
which has the energy function 
\begin{equation}\label{energy2}
E = \frac{1}{2}(\dot v)^2 - \frac{(n-2)^2}{2} v^2 
+ \frac{(n-2) \Lambda}{2n} v^{\frac{2n}{n-2}}.
\end{equation}
Solutions will lie on level sets of $E$ in the phase plane; 
see the discussion in \cite{Sc} for more details. We sketch 
some level sets of the energy function (\ref{energy2})
below. 

\begin {figure}[h]
\begin {center}
\includegraphics [width=4in]{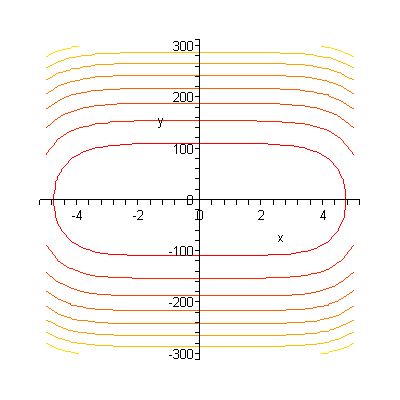}
\label{ball-level-sets}
\caption{Level sets for the energy function
(\ref{energy2}) with $n=3$ and $p=6$.}
\end {center}
\end {figure}

\section {Comparisons}
\label {comparison-sec} 

In this section we prove some basic 
comparison principles for minimizers 
of $\Phi_p$. The first such comparison 
is domain monotonicity. 
\begin {prop} 
If $D_1 \subset D_2\subset \R^n$ are bounded 
domains and $p\geq 1$ then $\calC_p(D_1) \geq
\calC_p(D_2)$.
\end {prop}

\begin {proof} 
This follows from the fact that $L^p(D_1) 
\cap W^{1,2}_0(D_1) \subset L^p(D_2) 
\cap W^{1,2}_0(D_2)$. 
\end {proof} 

Next we fix the domain $D$ and vary $p$. 
\begin {prop} \label {vary-p}
Fix a smooth, bounded domain $D\subset \R^n$, and let $\calC_p(D)$ 
be given by \eqref{C-p} for $p \geq 1$.  
Then the function $p \mapsto \calC_p(D)$ is continuous. 
\end {prop}

\begin {proof} Again, we use the fact that 
$C^\infty_0(D)$ is dense in $L^p(D) \cap 
W^{1,2}_0(D)$ and take $u$ to be a smooth 
function. In this case, the function 
$$p \mapsto \left ( \int_D u^p dV \right 
)^{2/p}$$
is smooth, so $\Phi_p(u)$ is a smooth 
function of $p$ for fixed $u \in C^\infty_0(D)$. 
The proposition follows from the definition 
of $\calC_p$. 
\end {proof} 

Our main result in this section is Theorem \ref{holder-comparison}, 
which reads 
$$
V^{2/p}\calC_p(D) > V^{2/q}\lambda_q(D)
$$ 
for $1\leq p<q$, where $V = \vol(D)$. 
\begin {proof} 
Recall that $C^\infty_0(D)$ is dense in both $L^p(D) \cap W^{1,2}_0(D)$ and 
$L^q(D) \cap W^{1,2}_0(D)$, so for the purposes of our comparison it will 
suffice to take $u \in C^\infty_0(D)$. In particular, $u \in L^p(D) \cap L^q(D)$. 
Use H\"older's inequality on the functions $u^p$ and $1$, with exponents 
$q/p$ and $q/(q-p)$, to obtain
$$
\left ( \int_D u^p dV\right )^{2/p} \leq \left [ V^{\frac{q-p}{q}} 
\left ( \int_D (u^p)^{q/p} dV \right )^{p/q} \right ]^{2/p}
= V^{\frac{2(q-p)}{qp}} \left ( \int_D u^q dV \right )^{2/q}.
$$
Then, by the variational character of $\calC_p$ and $\calC_q$, we have
$$
\calC_p(D) 
= \inf \left \{ 
\frac{\int_D |\nabla u|^2 dV}{\left ( \int_D u^p dV \right )^{2/p}} \right \} 
\geq V^{\frac{2(p-q)}{pq}} \inf \left \{ \frac{\int_D |\nabla u|^2 dV }{\left 
( \int_D u^q dV \right )^{2/q} }\right \}
= V^{\frac{2(p-q)}{pq}} \calC_q(D),
$$
which gives the desired inequality. Finally, we can only have equality in
H\"older's inequality if $u$ is constant, which (by the boundary conditions)
would force $u$ to be identically zero. 
This is impossible, and the inequality above must be strict. 
\end {proof} 

In dimension two one may take the limit as $p\to \infty$, to obtain 
$$
\calC_\infty (D) = \inf \left \{ \frac{\int_D |\nabla u|^2 dV}
{\|u\|_{L^\infty(D)}^2}: u \in L^\infty(D) \cap W^{1,2}_0(D), 
u \not \equiv 0 \right \} = \lim_{p \rightarrow \infty} \calC_p(D).
$$
By the monotonicty of $V^{2/p}\calC_p(D)$, this limit exists and is finite. 
Taking the limit as $p\to\infty$ in the scaling law \eqref{scaling-law2} with
$n=2$ shows that $\calC_\infty(rD) = \calC_\infty(D)$ for each positive $r$. 
For a fixed domain $D$, one can find disks $r\D$ and $R\D$ such that 
$r\D \subseteq D \subseteq R\D$, so that 
$\calC_\infty(R\D) \leq \calC_\infty(D) \leq \calC_\infty(r\D)$ 
by domain monotonicity. 
On the other hand, $\calC_\infty(r\D) = \calC_\infty(R\D)$ by scale invariance. 
We see that $\calC_\infty(D)$ does not depend on the domain at all, and write
$\calC_\infty$ for its common value for all planar domains. 

We close this section with a preliminary estimate for $\calC_\infty$.  
\begin {prop} \label{q=infty}
 $\calC_\infty \leq \pi$. 
\end {prop} 

\begin {proof} 
We may take $D$  to be the disk of radius $1$ centered at $0$. 
For any positive $\delta<1$, define the radial function 
$$
u(r) = \left \{ \begin {array}{cc} 1, & r<\delta \\ 
\frac{1 - r}{1 - \delta}, & 
\delta \leq r \leq 1
\end {array} \right. 
$$
Then $\| u\|_{L^\infty(D)} = 1$. Moreover,  
$$
\int_D |\nabla u|^2 dV = 
2\pi \int_\delta^1 \left( \frac{1}{\delta - 1} \right )^2 r\,dr = 
\pi \frac{1 + \delta}{1-\delta}.
$$
so that 
\begin{equation*}
\lim_{\delta \rightarrow 0^+} \left [ 
\int_D |\nabla u|^2 dV \right ] 
= \pi. \tag*{\qedhere}
\end{equation*}
\end {proof}

\section {Extremal domains} 
\label {extreme-sec} 

In this section we characterize the domains which are maxima or minima for 
$\calC_p$ under various constraints. We begin with a proof of 
Theorem \ref{faber-krahn}, that the ball uniquely minimizes $\calC_p(D)$ 
among all domains with a fixed volume. The proof follows the standard
proof of the Faber-Krahn inequality by symmetrization. 

\begin {proof} 
Let $u$ be a test function for $\Phi_p$. Without loss of generality, we can take 
$u \in C^\infty_0(D)$ and let 
$$
m = \min_{x \in D} \{u(x)\}, \qquad M = \max_{x \in D} \{u(x)\}.
$$
For $m \leq t \leq M$ let $D_t = \{ u > t \}$. 

Now we define a comparison function $u_*:B\rightarrow [m,M]$
as follows. First let $B_t$ be the ball centered at the origin with 
$\vol(B_t) = \vol(D_t)$. Then let $u_*$ be the radially symmetric function 
such that $B_t = \{ u_* > t\}$. By the co-area formula, 
$$
\int_{t}^{M} \int_{\del D_\tau} \frac 
{dA}{|\nabla u|} d\tau = \vol(D_t) = \vol(B_t) 
= \int_t^{M} \int_{\del B_\tau} \frac
{dA}{|\nabla u_*|} d\tau.
$$
Differentiating with respect to $t$ gives  
\begin {equation} \label{eq-deriv}
\int_{\del D_t} \frac{dA}{|\nabla u|} = 
\int_{\del B_t} \frac{dA}{|\nabla u_*|}
\end {equation}
for all $t$. Then 
\begin {eqnarray} \label{fubini}
\int_D u^p dV & = & \int_m^M \int_{\del D_t}
\frac{u^p dA}{|\nabla u|} dt = \int_m^M
t^p \int_{\del D_t} \frac{dA}{|\nabla u|}  dt
\nonumber \\
& = & \int_m^M t^p \int_{\del B_t}\frac 
{dA}{|\nabla u_*|} dt = \int_B u_*^p dV.
\end {eqnarray}
Now, for $m \leq t \leq M$ let 
$$
\psi(t) = \int_{D_t} |\nabla u|^2 dV, \qquad 
\psi_*(t) = \int_{B_t} |\nabla u_*|^2 dV.
$$
By the co-area formula 
$$
\psi' = -\int_{\del D_t} |\nabla u| dA, \qquad 
\psi_*' = -\int_{\del B_t} |\nabla u_*| dA.
$$
We use the Cauchy-Schwarz inequality, the isoperimetric inequality, 
and the fact that the normal derivative of $u_*$ is constant on 
$\del B_t$ to see
\begin {eqnarray*}
\left ( \int_{\del D_t} |\nabla u| \,dA \right ) \left ( \int_{\del D_t} \frac{dA}
{|\nabla u_*|} \right ) 
& \geq & \left (\int_{\del D_t} \,dA \right )^2 = \left( \Area(\del D_t)\right )^2 \\
& \geq & \left ( \Area (\del B_t) \right )^2
= \left ( \int_{\del B_t} |\nabla u_*| \,dA 
\right ) \left ( \int_{\del B_t} \frac{dA}{|\nabla u_*|} \right ).
\end {eqnarray*}
We use equation (\ref{eq-deriv}) to cancel the common factor of 
$$
\int_{\del D_t} \frac{dA}{|\nabla u|} = \int_{\del B_t} \frac{dA}{|\nabla u_*|},
$$
and so 
$$
-\psi' = \int_{\del D_t} |\nabla u| \,dA \geq \int_{\del B_t} |\nabla u_*| \,dA 
= -\psi_*'.
$$
Integrating this last differential inequality and using 
$\psi(M) = 0 = \psi_*(M)$ we see that 
$$
\int_D |\nabla u|^2 dV = \psi (0) \geq \psi_* (0) = \int_B |\nabla u_*|^2 dV.
$$
This inequality, combined with (\ref{fubini}) and (\ref{C-p}), give the desired 
inequality on the eigenvalues: 
$$
\calC_p(D) \geq \calC_p(B).
$$
Moreover, equality of the eigenvalues forces all level sets $\del D_t$ to 
be spheres centered at the origin. Also, the equality case of the 
Cauchy-Schwarz inequality forces $|\nabla u|$ to be constant on the 
level set $\del D_t$. Thus $u$ must be radially symmetric,  
and so in this case $u = u_*$. 
\end {proof} 

Next we fix the inradius $R(D)$ of the domain rather than the volume, 
where $R(D)$ is the supremum radius of all balls contained in $D$. 
\begin {lemma}
Among all bounded domains $D \subset 
\R^n$ with a fixed inradius, the ball 
maximizes $\calC_p$ for all $p \geq 1$. 
\end {lemma} 
\begin {proof} 
If the inradius of $D$ is $R$, then $D$ contains a ball of radius $r$ for each $r<R$.
 The result now follows from domain monotonicity. 
\end {proof}  

\begin {prop} 
Let $p\geq 1$.  Let $D$ be a  smooth, convex domain $D \subset \R^n$ with inradius $R$. Let $u>0$ solve (\ref{p-pde}) on $D$ and let 
$u_M = \max \{u(x) : x \in D\}$. Then 
\begin {equation} \label{u-max-inradius}
u_M^{2-p} \leq \frac{2\Lambda R^2} {p A_p^2}, 
\end {equation}
where 
$$
A_p = \int_0^1 \frac{dt}{\sqrt{1-t^p}}
= \sqrt{\pi} \frac{\Gamma(1+\frac{1}{p})}{\Gamma(\frac{1}{2} + \frac{1}{p})}.
$$
Moreover, (\ref{u-max-inradius}) is an equality in the case of a slab.
\end {prop} 

\begin{rmk} By the discussion in the introduction, 
the maximum $u_M$ is well-defined. \end {rmk} 

Hersch \cite[Th\'eor\`eme~8.1]{Her} proved $\lambda(D) \geq \pi^2/(4R^2)$ 
for the fundamental frequency while Sperb \cite{Sp}  proved $u_M\leq R^2$ 
for the maximum value of the torsion function, in each case 
for a convex domain of  inradius $R$. 
The estimate \eqref{u-max-inradius} is a common generalisation of these
results. More refined results are known in the cases $p=1$ and $p=2$
 (see \cite{MH}, for example). 

\begin {proof} We follow Section 6.2.2 of \cite{Sp}. The $P$-function
$$
v(x) = |\nabla u(x)|^2 + \frac{2\Lambda}{p} u^p(x),
$$
introduced by Payne, assumes its maximum at the point where $u$ assumes 
its maximum. Thus 
$$
|\nabla u(x)|^2 + \frac{2\Lambda}{p}
u^p(x) \leq \frac{2\Lambda}{p} u_M^p,
$$
which we can rearrange to read 
\begin{equation} \label{u-max-de} 
|\nabla u(x)| \leq \sqrt{\frac{2\Lambda}{p}} \sqrt{u_M^p - u^p(x) }. 
\end{equation}
Let $\delta_D(P)$ be the distance from the point $P$ where $u$ assumes its 
maximum to the boundary of $D$ and integrate (\ref{u-max-de}) 
along a line segment which starts at $P$ and terminates  
at a point on $\del D$ closest to $P$. Then 
$$
R(D) \geq \delta_D(P) \geq \sqrt{\frac{p}{2\Lambda}} u^{(2-p)/2} A_p
\quad \text{where } A_p 
= \int_0^1 \frac{dt}{\sqrt{1-t^p}}.
$$
The inequality (\ref{u-max-inradius}) follows. Moreover, if $D$ is a 
slab then (\ref{slab-1st-integral}) implies (\ref{u-max-de}) is actually 
an equality, and so (\ref{u-max-inradius}) is also an equality. 
\end {proof}

\section {Open questions}
\label {open-sec}

In this final section we collect a small sample of interesting, related questions. 

In the present paper, we have restricted our attention to the minimizer of 
the the functional $\Phi_p$, which corresponds to the bottom of the 
spectrum of the eigenvalue equation 
$$
\Delta u + \lambda u^{p-1} = 0.
$$
This functional should have other critical points above the minimizer, 
which correspond to sign changing solutions of the boundary value 
problem (\ref{p-pde}). What can one say about these higher order eigenvalues? 
If $n\geq 3$ and $p< \frac{2n}{n-2}$ (or if $n=2$) 
is the spectrum discrete? Is there a sequence 
of eigenvalues $\lambda_{p,j}$, with $j=1,2,3,\dots$, 
such that 
$$0 < \lambda_{p,1} < \lambda_{p,2} \leq 
\lambda_{p,3} \leq \cdots \rightarrow \infty ?$$

We return to the minimum $\calC_p(D) = \inf \{\Phi_p(u)\}$. 
In dimension $n \geq 3$, does the limit of $\calC_p(D)$ exist as 
$p \rightarrow \frac{2n}{n-2}$ from below? It would be nice to 
characterize the domains on which the limit exists, in terms of geometry. 
We showed in Proposition \ref{vary-p} that the eigenvalue $\calC_p$ is a 
continuous function of $p$. Is the same true of the eigenfunction? 
We have shown that, in dimension two, $\lim_{p \to \infty} \calC_p(D) 
= \calC_\infty$ exists and is independent of the domain. Can the bound 
$\calC_\infty\leq \pi$ be improved or is it sharp? We proved the inequality 
(\ref{u-max-inradius}) is realized in the case of a slab. Are there any other
domains for which (\ref{u-max-inradius}) is an equality?

\end{document}